\documentclass[letterpaper, 12pt]{article}

\usepackage[letterpaper, margin=1in, bmargin=1.1in]{geometry}  % set the margins to 1in on all sides
\usepackage{amsmath}               % great math stuff
\usepackage{amsfonts}              % for blackboard bold, etc
\usepackage{amsthm}                % better theorem environments
\usepackage[colorlinks]{hyperref}
\usepackage{natbib}
\usepackage{enumerate}
\newcommand{\numalph}{\renewcommand{\labelenumi}{\hspace*{1em}\alph{enumi})\quad}}

\theoremstyle{definition}

\newtheorem{theorem}{Theorem}[section]
\newtheorem{lemma}[theorem]{Lemma}

\newtheorem{cor}[theorem]{Corollary}
\newtheorem{defin}[theorem]{Definition}

\newtheorem{rem}[theorem]{Remark}
\newtheorem{rems}[theorem]{Remarks}

\renewcommand{\v}[1]{{\mathbf{#1}}}
\newcommand{\inner}[2]{\langle #1,#2\rangle}

\newcommand{\tren}{\hspace{0pt}}
\newcommand{\R}{\mathbb{R}}

\newcommand{\phhi}{\varphi}
\renewcommand{\phi}{\varphi}

\newcommand{\Kal}[1]{{\mathcal{#1}}}

\newcommand{\keywords}[1]{\par\noindent{\small{\em Keywords\/}: #1}}

\renewcommand{\em}{\sffamily\slshape}

\begin{document}

\title{Frenet Curves and Successor Curves: Generic Parametrizations of the Helix and Slant Helix}

\author{Toni Menninger\\\normalsize E-mail: toni.menninger@gmail.com}

%\thanks{E-mail: toni.menninger@gmail.com}}

\maketitle

\begin{abstract}
%\noindent
In classical curve theory, the geometry of a curve in three dimensions is essentially characterized by their invariants, curvature and torsion. When they are given, the problem of finding a corresponding curve is known as 'solving natural equations'. Explicit solutions are known only for a handful of curve classes, including notably the plane curves and general helices. 

This paper shows constructively how to solve the natural equations explicitly for an infinite series of curve classes. For every Frenet curve, a family of successor curves can be constructed which have the tangent of the original curve as principal normal. Helices are exactly the successor curves of plane curves and applying the successor transformation to helices leads to slant helices, a class of curves that has received considerable attention in recent years as a natural extension of the concept of general helices. 

The present paper gives for the first time a generic characterization of the slant helix in three-dimensional Euclidian space in terms of its curvature and torsion, and derives an explicit arc-length parametrization of its tangent vector. These results expand on and put into perspective earlier work on Salkowski curves and curves of constant precession, both of which are subclasses of the slant helix.  

The paper also, for the benefit of novices and teachers, provides a novel and generalized presentation of the theory of Frenet curves, which is not restricted to curves with positive curvature. Bishop frames are examined along with Frenet frames and Darboux frames as a useful tool in the theory of space curves. The closed curve problem receives attention as well.

%\vspace{4 mm}
\end{abstract}

\keywords{Curve theory; General helix; Slant helix; Curves of constant precession; Frenet equations; Frenet curves; Frenet frame; Bishop frame; Darboux frame; Natural equations; Successor curves; Salkowski curves; Closed curves.}

\newpage

\section{Introduction}

In  classical three-dimensional curve theory, the geometry of a curve is essentially characterized by two scalar functions, {\em curvature} $\kappa$ and {\em torsion} $\tau$, which represent the rate of change of the tangent vector and the osculating plane, respectively. Given two continuous functions of one parameter, there is a space curve (unique up to rigid motion) for which the two functions are its curvature and torsion (parametrized by arc-length). The problem of finding this curve is known as {\em solving natural} or {\em intrinsic equations} and requires solving the Frenet (or an equivalent) system of differential equations (see \cite{scofield}; \cite{Hoppe1862}). One actually solves for the unit tangent and obtains the position vector by integration. Explicit solutions are known only for a handful of curve classes, including notably the plane curves (solved by Euler in 1736)
%\footnote{cf. Eric W. Weisstein, Natural Equation. From MathWorld--A Wolfram Web Resource. http://mathworld.wolfram.com/NaturalEquation.html})
and general helices. 

Recall that {\em general helices} or {\em slope lines} are defined by the property that their tangent makes a constant angle with a fixed direction in every point. Their centrode, the momentary axis of motion, is fixed. Similarly, the principal normal vector of a {\em slant helix} makes a constant angle with a fixed direction and its centrode precesses about a fixed axis. The term was only recently coined (\cite{izumya}) but such curves have appeared in the literature much earlier (e. g. \cite{Hoppe1862}; \cite{bilinski1955}: 180; \cite{Hoschek1969}). In particular, Salkowski studied slant helices with constant curvature (\cite{salkowski1909}, \cite{monterde}) and Scofield derived closed-form arc-length parametrizations for {\em curves of constant precession}, slant helices with the speed of precession (and not just the angle) being constant (\cite{scofield}). 
\cite{menninger} (chapter 9) appears to have been the first study of what is now known as slant helices as a curve class in its own right. More recently, there has been considerable interest in the slant helix among curve geometers (e.g. \cite{kula2005}, \cite{kula2010}, \cite{ali}, \cite{Camci2013}).

The close relation between slant helix and general helix leads to the concept of the successor transformation of Frenet curves. Given a Frenet moving frame, we can construct a new Frenet frame in which the original tangent vector plays the role of principal normal. It turns out that general helices are the successor curves of plane curves, and slant helices are the successor curves of general helices. Curves of constant precession are the successor curves of circular helices, which in turn are the successor curves of plane circles. 
This transformation could be repeated indefinitely to give rise to new, yet unexplored classes of curves. 
 
In this paper, after a succinct review of Frenet curve theory, the successor transformation of Frenet curves is introduced and its general form given (section \ref{sectionsuccessorcurves}). It is then, in section \ref{sectionapplied}, applied to plane curves and helices to derive a generic characterization of the slant helix.

\newpage
\section{Frenet Space Curves}
\subsection{Regular Curves and Moving Frames}
This section introduces basic concepts of classical curve theory. The following will be concerned with regular curves in oriented three-dimensional euclidian space  $\R^3$, endowed with an {\em inner product},
\[\inner{\v u}{\v v} =\|\v u\|\|\v v\|\cos \theta,\] and a {\em cross product}, \[\v u \times \v v=\|\v{u}\|\|\v v\|\sin \theta \,\v N,\] for any $\v{u}, \v v \in \R^3$, where $\theta \in [0, \pi]$ is the angle between $\v{u}$ and $\v v$ and $\v N$ is the positively oriented unit vector perpendicular to the plane spanned by $\v{u}$ and $\v v$. In particular, \[\v{u}\perp \v v \Leftrightarrow \inner{\v{u}}{\v v}=0 \ \text{and} \ \v{u} \parallel \v v \Leftrightarrow \v{u} \times \v{v}=0.\]

\begin{defin}[Regular Curves]

A  {\em parametrized curve} is a $C^{k}$ map \,$\v{x}:I\mapsto \R^n$ (where $k \geq 0$ denotes the differentiation order) defined on an interval $I\subset\R$. We can picture it as the trajectory of a particle moving in space, with position vector $\v{x}(t)$ and (if $k\geq 1$) velocity $\dot{\v{x}}(t)=\frac{dx}{dt}$ at time $t$. The image $\v{x}[I]$ is called {\em trace} of the parametrized curve.
 
A parametrized $C^1$ curve $\v{x}(t)$\, is called {\em regular} \,if it has a nonvanishing derivative $\dot{\v{x}}(t)$ at every point. A regular parametrized curve $\v{x}(t)$ has a well-defined and continuous unit tangent vector $\v T(t)=\dot{\v{x}}(t)/ \|\dot{\v{x}}(t)\|$ and can be parametrized by arclength (cf. \cite{kuehnel}: 9), denoted $s$ throughout this paper. The intrinsic geometric properties of regular curves are independent of regular parameter transformations and any regular curve can be represented by an arclength (or {\em unit speed}) parametrization $\v{x}(s)$ with $\v{x}'=\v T$. 
\end{defin}

\begin{defin}[Vector Field]
A $C^k$ {\em vector field} along a parametrized curve $\v{x}:I \subset \R\mapsto \R^n$ is a $C^k$ map $\v{v}: I \mapsto \R^n$, and a {\em unit vector field} is a vector field of unit length, $\v{V}: I \mapsto S^{n-1}$. A vector field $\v{v}$ that is everywhere linearly dependent on a unit vector field $\v{V}$ of the same differentiation order is called {\em normalizable}.
It is $\v{v(s)}=\inner{\v{v(s)}}{\v{V(s)}}\cdot \v{V(s)}$ and it is said that $\v V$ normalizes $\v v$.
Any nowhere-vanishing (non-degenerate) vector field is normalizable but that is not in general the case. 
A vector field $\v v$ is called {\em tangential}\, if $\,\v v(s) \parallel \v T(s)$ and {\em normal} \,if  \, $\v v(s) \perp \v T(s)$ for all $s \in I$. A normal $C^1$ vector field is called {\em parallel} to the curve if its derivative is tangential along the curve. 
\end{defin}

\begin{lemma}
\label{skewsymm}
For any two differentiable unit vector fields $\v{V}$ and $\v{W}$ keeping a constant angle, it is 
\[\inner{\v{V}}{\v{W}}=const. \;\Leftrightarrow \;\inner{\v{V}}{\v{W}}'=0 \;\Leftrightarrow \; \inner{\v{V}'}{\v{W}}= -\inner{\v{V}}{\v{W}'}.\] 
Any $C^1$ unit vector field is orthogonal to its derivative. In particular, $\v T'$ is normal.
\end{lemma}

\begin{defin}[Canonical Curvature and Principal Normal]
Let $\v{x}=\v{x}(s)$ be a unit speed \,$C^{2}$-{\tren}curve and \,$\v T=\v{x}'$\, its unit tangent vector. Then the continuous function \[ \kappa_+ =\|\v T'\|=\|\v{x}''\|\] is called its  {\em canonical curvature}. A point with $\kappa_{+}(s)=0$ is called {\em inflection point}, otherwise {\em strongly regular}. A  curve without inflection points is also called  {\em strongly regular}. 

The trace of a (sub)interval of the arc length parameter is called a {\em curve segment}. A {\em (straight) line segment} is a curve segment with vanishing curvature. 

In strongly regular points, the {\em canonical principal normal vector} is defined as \[\v{N}_{+} =\kappa_{+}^{-1}\;\v{T}'.\]
\end{defin}

\begin{rem}Every regular point of a space curve in $\R^3$ has a well-defined {\em tangent line}, spanned by $\v{T}$, and perpendicular  {\em normal plane}. 
Every strongly regular point has a well-defined {\em principal normal line}  spanned by $\v{N}_{+}$, {\em osculating plane } spanned by $\v{T}$ and $\v{N}_{+}$, and {\em rectifying plane} perpendicular $\v{N}_{+}$. 
\end{rem}

\begin{defin}[Tangential Moving Frame]
A {\em moving frame}  along a regular curve is a tuple $(\v V_1 ... \v V_n)$ \,of $n$ differentiable, pairwise orthogonal unit vector fields forming a {\em positively oriented, orthonormal basis} of $\R^n$ attached to every point of the curve. It is called {\em tangential} if $\v V_1 = \v T$. Two tangential moving frames are called {\em equivalent} if their tangential components are the same.
\end{defin}

\noindent {\em The remainder of this paper is constricted to unit speed curves and moving frames in $\R^3$.}

\begin{theorem}[Differential Equations for Moving Frame and Fundamental Theorem]
\label{diffeq}
Let $(\v T, \v N, \v M)$ be a tangential moving frame of a regular space curve. Then a system of differential equations, with coefficient functions $k_1, k_2, k_3:I\mapsto \R$ holds:
\[ \begin{pmatrix}\v{T}\\ \v{N}\\ \v{M} \end{pmatrix}'
=\begin{pmatrix} \; 0&\;k_1&k_2 \\ -k_1&\;0&k_3 \\ -k_2 &-k_3&\;\,0 \end{pmatrix}\cdot
\begin{pmatrix} \v{T}\\ \v{N}\\ \v{M} \end{pmatrix}. \]

Any three continuous functions $k_1(s), k_2(s), k_3(s):I \subset \R\mapsto\R$ determine a unique (up to a rotation) moving frame of which they are the differential equation coefficients, which in turn is tangential to a uniquely determined (up to a translation) regular $C^2$ curve. In consequence, three arbitrary continuous functions taken as coefficient functions of a tangential moving frame determine a regular unit speed space curve up to a rigid Euclidian motion.
\end{theorem}

\begin{proof} 
Lemma \ref{skewsymm} implies that differentiating a moving frame gives rise to a skew symmetric matrix of coefficient functions. The resulting system of linear differential equations is known to have a solution, which is a moving frame, and it is unique up to the choice of initial conditions (cf. \cite{kuehnel}: 28). Two moving frames with the same differential coefficient functions differ only by a rotation. Integrating the tangent component gives an arclength parametrization $$x(s)=x(s_0) +\int_{s_0}^s \v T(u)du$$ of a regular $C^2$ curve, unique up to a translation. Therefore, two solution curves are identical up to a rigid motion. 
\end{proof}

\begin{theorem}[Totality of Tangential Moving Frames]
\label{totality}
Let $(\v T, \v N, \v M)$ be a tangential moving frame of a regular space curve. The totality of equivalent tangential moving frames $(\v T, \overline {\v N}, \overline {\v M})$ can be expressed as follows, depending on an arbitrary differentiable function $\phhi$:
\[\begin{pmatrix} \v T \\ \v {\overline N}\\ \v {\overline M}\end{pmatrix}=
 \begin{pmatrix} 1&0&0 \\ 0& \cos\phhi& -\sin \phhi  \\0&  \sin\phhi &  \quad\cos \phhi\end{pmatrix}\cdot
\begin{pmatrix} \v T \\ \v {N}\\ \v { M}\end{pmatrix}\]
The differential equation coefficients get transformed as follows:
\[ \overline k_1=k_1 \cos\phhi - k_2 \sin \phhi, \quad \overline k_2=k_1 \sin\phhi + k_2 \cos \phhi, \quad \overline k_3=k_3-\phhi'.\]

\end{theorem}
\begin{proof}
The rotation matrix, denoted $R_1(\phhi)$, represents the subgroup of the special orthogonal group $SO(3)$ that leaves the first (tangential) component of the moving frame unchanged. It is obvious that all eqivalent tangential moving frames can be expressed in this way. Denoting $F=(\v T, \v N, \v B)^t$ and the coefficient matrix as $K(F)$, we have ${\overline F=R_1(\phhi)F}$, ${F'=K(F)F}$ and \,$\overline F'=R_1'(\phhi)F+R_1(\phhi)F'=
[R_1'(\phhi)+R_1(\phhi)K(F)]R_1(-\phhi)\,\overline F=K(\overline F)\overline F$. Evaluation of $K(\overline F)$ confirms the stated relationships.
\end{proof}

\subsection{A Generalization of the Theory of Frenet Space Curves}
Conventionally, the theory of Frenet curves in three dimensions is restricted to strongly regular $C^3$ curves and introduced by constructing the Frenet moving frame \,$(\v{T}, \v{N}_+, \v{T} \times \v{N}_+)$ (e.~g. \cite{kuehnel}: 13). 
In this section, a generalized version of Frenet theory is presented which allows that important classes of curves, including helices and slant helices (section 3 of this paper), can be treated as Frenet curves even if they have inflection points or even line segments. The price to be paid for relaxing the assumption of strong regularity is the loss of uniqueness of the Frenet invariants.

\begin{defin}[Frenet Curves and Frenet Apparatus]
\label{frenetcurves}

Given a unit speed $C^{2}$ space curve $\v{x}:I\subset \R\mapsto \R^3$ with unit tangent vector ${\v{T}:I\mapsto S^2}$,\; ${\v{T}=\v{x}'}$. 
\begin{enumerate}\numalph

\item A {\em principal normal} is a normal $C^1$ unit vector field  $\v{N}:I\mapsto S^2$ normalizing $\v T'$: \\ $\v N \perp \v T$ and \,$\v N \parallel \v T'$.

\item A {\em binormal} is a normal $C^1$ unit vector field $\v{B}:I\mapsto S^2$ with a normal derivative:\, $\v B \perp \v T$ and \,$\v B' \perp \v T$ (which implies  \,$\v B \perp \v T'$).

\item A {\em Frenet frame} of the curve is a tangential moving frame \,$(\v{T}, \v{N}, \v{B})$ where $\v{T}$ is its tangent, $\v{N}$ is a principal normal and $\v{B}$ is a binormal.

\item Given a Frenet frame  $(\v{T}, \v{N}, \v{B})$, the following are well-defined and continuous:
\begin{alignat*}{10}
	& &&\bullet\; \text{The (signed) {\em curvature} } && \kappa:&&I\mapsto \R, &&\kappa &&=\; \inner{\v{T}'}{\v{N}} \\
	& &&\bullet\; \text{The {\em torsion} } && \tau:&&I\mapsto \R, &&\tau&&=\; \inner{\v{N}'}{\v{B}}\\
	& &&\bullet\; \text{The {\em Lancret curvature} } && \omega_+:\;&&I\mapsto \R, \quad\quad&&\omega_+&&=\; \sqrt{\kappa^2+\tau^2}\\
	& &&\bullet\; \text{The {\em Darboux vector field} \quad} && \v D_\omega:&&I\mapsto \R^3, &&\v D_\omega&&=\; \tau \v{T}+\kappa \v{B}
\end{alignat*}

\item The tuple $\Kal F = (\v{T}, \v{N}, \v{B}, \kappa, \tau)$ is called a {\em Frenet apparatus} or  {\em Frenet system} and the pair $(\kappa, \tau)$ a {\em Frenet development} associated with the curve.

\item A Frenet apparatus or moving frame is called {\em generic} if along every line segment, curvature and torsion both vanish (or if it has no line segments). 

\item A curve is called a {\em Frenet curve} if it has a Frenet frame, and a {\em generic Frenet curve} if it has a generic Frenet frame.
\end{enumerate}
\end{defin}

\begin{rem}[Existence of Frenet Apparatus]
\cite{wintner} defined a Frenet curve as one that possesses a binormal, while \cite{nomizu} and \cite{wonglai} defined it equivalently as a curve possessing a Frenet moving frame.
The existence of either a principal normal or a binormal is sufficient for the existence of a Frenet frame since the third component can be constructed as  $\v{B}=\v{T} \times \v{N}$ or $\v{N}=\v{B} \times \v{T}$, respectively.
A regular curve is a Frenet curve if and only if $\v{T}'$ is normalizable, or equivalently, if the map from each strongly regular point to the osculating plane has a continuous extension. Curves that discontinuously  "jump" between planes cannot be Frenet-framed. 
By definition, a Frenet curve is regular and at least of order $C^2$. The components of the Frenet frame are at least of order  $C^1$ and the Frenet development at least of order~$C^0$. 
\end{rem}

\begin{rem}[Surface Curves and Darboux Moving Frame]
\label{darbouxframe}
A special case of a tangential moving frame is the {\em Darboux Moving Frame} $(\v T, \v S, \v {\hat N})$ of a $C^2$ curve with tangent $\v T$, situated on a regular parametrized $C^2$ surface with unit surface normal $\v {\hat N}$, and $\v S=\v {\hat N}\times \v T$. The differential equation coefficients of the Darboux frame are known as {\em geodesic curvature} $\kappa_g=\inner{\v T'}{ \v S}$,  {\em normal curvature} $\kappa_n=\inner{\v T'}{ \v {\hat N}}$ and {\em geodesic torsion} $\tau_g=\inner{\v S'}{ \v {\hat N}}$. 

{\em Geodesics} are surface curves with $\kappa_g\equiv 0$. They are Frenet curves because the unit surface normal is a principal normal. Curves with $\kappa_n\equiv 0$ are  {\em Asymptotic curves} and their Darboux frame is a Frenet frame. Finally, curves with $\tau_g\equiv 0$ are called {\em Lines of curvature}.
\end{rem}

The Frenet moving frame is the special case of a tangential moving frame with vanishing coefficient $k_2$. The Frenet differential equations and the Fundamental Theorem for Frenet curves follow directly from definition \ref{frenetcurves} and theorem \ref{diffeq}.

\begin{theorem}[Frenet Equations]
\label{freneteq}
Any Frenet apparatus $(\v{T}, \v{N}, \v{B}, \kappa, \tau)$ satisfies the {\em Frenet equations}
\[\v{T}'=\kappa \v{N}, \qquad \v{N}'=-\kappa \v{T} + \tau \v{B},  \qquad \v{B}'=-\tau \v{N}.\]
In matrix notation, we have
\[ \begin{pmatrix}\v{T}\\ \v{N}\\ \v{B} \end{pmatrix}'
=\begin{pmatrix} \; 0&\;\kappa&\;\,0 \\ -\kappa&\;0&\;\,\tau \\ \;0 &-\tau&\;\,0 \end{pmatrix}\cdot
\begin{pmatrix} \v{T}\\ \v{N}\\ \v{B} \end{pmatrix}=
\v D_\omega \times \begin{pmatrix}\v{T}\\ \v{N}\\ \v{B}\\\end{pmatrix}. \]
We also have the relationships
\[\kappa^2=\kappa_+^2=\inner{\v{T}'}{\v{T}'}, \qquad \omega_+^2=\inner{\v{N}'}{\v{N}'},\qquad \tau^2=\inner{\v{B}'}{\v{B}'}.\]
\end{theorem}

\begin{rem}[Darboux Rotation Vector]
The right hand side of the Frenet equation suggests that the {\em Darboux rotation vector} \,$\v D_\omega$ can be interpreted as the angular momentum vector of the Frenet frame. Its direction determines the moving frame's momentary axis of motion (its {\em centrode}) and its length the angular speed $\|D_\omega\|=\omega_+$. When $D_\omega$ is normalizable, there is a continuous {\em unit Darboux vector} $D$ satisfying \,$D_\omega=\omega D_\omega$, with \,$|\omega|=\omega_+$. The function $\omega$ is the {\em signed Lancret curvature}.
\end{rem} 

\begin{theorem}[Fundamental Theorem for Frenet Curves]
\label{fundamental}
Given a pair of continuous functions $\kappa(s)$ and $\tau(s): I \subset \R\mapsto\R$. Then there exists a unit speed Frenet curve with $\kappa(s)$ and $\tau(s)$ as curvature and torsion. Two Frenet curves with the same Frenet development are congruent and differ only by a rigid Euclidian motion.
The equations $\kappa=\kappa(s)$ and $\tau=\tau(s)$ are also know as the {\em natural} or intrinsic equations of the associated space curve. 
\end{theorem} 

\begin{theorem}[Equivalence of Frenet Developments]
Two Frenet developments $(\kappa, \tau)$ and $(\overline{\kappa}, \overline{\tau})$ characterize the same Frenet curve if and only if there exists a $C^1$-Funktion $\phhi$ so that
\[\kappa\cos\phhi=\overline\kappa, \quad \kappa\sin\phhi=0    \quad
\text{and} \quad \phhi'=\tau-\overline\tau.\]  
\end{theorem}
\begin{proof}
The conditions (due to \cite{wonglai}: 11) can be read off theorem \ref{totality} by setting $k_1=\kappa, k_2=0, k_3=\tau$.
\end{proof}

\begin{cor}[Curves without Line Segments]
For a Frenet curve without line segments (i.~e. $\kappa_+>0$ almost everywhere), \:$ \phhi=k \pi=const.$ \,due to continuity. It has exactly two equivalent Frenet apparatuses \: $(\v{T}, \v{N}, \v{B}, \kappa, \tau)$ \; and \: $(\v{T}, -\v{N}, -\v{B}, -\kappa, +\tau)$. Its curvature is unique up to sign and it has a unique and well-defined torsion. 
\end{cor}

\begin{cor}[Total Torsion] 
\label{totaltorsion}
The {\em total torsion} of a curve or curve segment, \, $\int_{0}^l\tau(s)\, ds$, 
is invariant modulo $\pi$ \,{\em if} \,none of its end points lies within a line segment, because \:$\int_{0}^l\tau-\int_{0}^l\overline\tau=\phhi(l)-\phhi(0)$ \;and \:$\,\sin\phhi(l)=\sin\phhi(0)=0$. 
\end{cor}

\begin{rem}[Hierarchy of Frenet Curves]
It is useful to distinguish the following classes of curves:
\begin{description}

\item[Strongly regular curves:] Every strongly regular $C^{3}$-{\tren}curve is a Frenet curve. It has a canonical Frenet apparatus with principal normal $\v{N}=\v{N}_{+}$ and curvature $\kappa=\kappa_{+}>0$. 

\item[Analytic curves:] \cite{nomizu} showed that an analytic curve is a Frenet curve and has at most isolated inflection points, unless it is a line segment. In fact, Nomizu showed constructively that this is the case for every {\em normal} curve, that is $C^\infty$-curve which in every point has at least one non-vanishing derivative.

\item[Frenet curves without line segments:] For a Frenet curve without line segments,  $\v{N}$, $\v{B}$, and $\kappa$ are unique up to sign and $\tau$ is a curve invariant. 

\item[Frenet curves with generic Frenet apparatus:] A generic Frenet apparatus is stationary along line segments (if such exist). It is unique up to the sign of $\v{N},\,\v{B}$, and $\kappa$. 

\item[Frenet curves:] A Frenet curve containing line segments has an infinite family of Frenet apparatuses. In strongly regular points, the Frenet frame is determined up to sign but along each line segment, it can rotate arbitrarily around the tangent and the torsion can vary arbitrarily. Only the magnitude of the curvature, \,$|\kappa|=\kappa_+$, is a curve invariant. The total torsion modulo $\pi$ over a curve segment is invariant (\ref{totaltorsion}).

\item[Regular curves:] A regular curve need not be a Frenet curve. However, any regular $C^2$-curve has a family of {\em Bishop frames} (\cite{bishop} and \cite{menninger}: ch. 6). 
\end{description}
\end{rem}

\subsection{The Bishop Frame}

\begin{defin}[Bishop Frame]
Given a unit speed $C^{2}$ space curve with unit tangent vector $\v T$.
A {\em Bishop frame} (or {\em parallel transport frame}) is a tangential moving frame  $(\v T, \v N_1, \v N_2)$ with both normal components being parallel, i.e. having tangential derivatives along the curve. Equivalently, it is a tangential moving frame with vanishing coefficient function $k_3$.
The equations for the derivatives assume the form 
%$$ \v T'=k_1 \v N_1 + k_2 \v N_2, \quad \v N_1'=-k_1 \v T, \quad \v N_2'=-k_2 \v T$$ 
\[ \begin{pmatrix} \v{T}\;\\ \v N_1\\ \v N_2\end{pmatrix}'
=\begin{pmatrix} 0&k_1&k_2 \\ -k_1&0&0 \\ -k_2 &0&0 \end{pmatrix}\cdot
\begin{pmatrix} \v{T}\;\\ \v N_1\\ \v N_2 \end{pmatrix}.\]
and the tuple $(k_1, k_2)$ is called the {\em Bishop (or normal) development} of the associated curve.
\end{defin} 

The Frenet and the Bishop frame are the only ways to construct a tangential moving frame with one coefficient of the derivative equations vanishing (cf. \cite{bishop}). A Frenet frame can be rearranged into a Bishop frame, and vice cersa.

\begin{theorem}[Rearrangement of Moving Frames]
\label{rearrange}
Given a Frenet system \, $(\v{T}, \v{N}, \v{B}, \kappa, \tau)$. 
\begin{enumerate}\numalph
\item
$(\v{T}, -\v{N}, -\v{B}, -\kappa, \tau)$ is an equivalent Frenet system. 

\item $(-\v{T}, -\v{N}, \v{B}, \kappa, -\tau)$ is also a Frenet system but it's not equivalent because orientation is not preserved.

\item $(\v{B}, -\v{N}, \v{T}, \tau, \kappa)$ is a Frenet system characterizing a curve with $\v B$ as unit tangent.

\item $(\v{N}, -\v{T}, \v{B}, \kappa, \tau)$ is a Bishop system with $\v N$ as unit tangent.

\item Conversely, a Bishop system  \,$(\v T, \v N_1, \v N_2, k_1, k_2)$ can be rearranged into a Frenet system 
 \,$(-\v N_1, \v T, \v N_2, k_1, k_2)$.
\end{enumerate}
\end{theorem}

\begin{theorem}[Family of Bishop Frames]
\label{bishopfamily}
Given a Bishop frame  \,$(\v T, \v N_1, \v N_2)$. The totality of equivalent Bishop frames consists of frames \,$(\v T, \v {\overline N}_1, \v{ \overline N}_2)$ of the form:
\[ \begin{pmatrix} \v {\overline N}_1\\ \v {\overline N}_2\end{pmatrix}=
 \begin{pmatrix} \cos\phhi_0& -\sin \phhi_0  \\ \sin\phhi_0 &  \quad\cos \phhi_0\end{pmatrix}\cdot
\begin{pmatrix} \v {N}_1\\ \v { N}_2\end{pmatrix}
\]
They differ only by a rotation by a constant angle $\phhi_0\in \R$ in the normal plane. 
%$$  \v {\overline N}_1=\quad\cos\phhi_0 \v N_1 + \sin \phhi_0 \v N_2$$ $$\v {\overline N}_2=-\sin\phhi_0 \v N_1 + \cos \phhi_0 \v N_2$$
The Bishop development of a regular curve, conceived of as a parametrized plane curve, is uniquely determined up to a rotation about the origin. 
\end{theorem}
\begin{proof}
This is a consequence of theorem \ref{totality} with $k_3=\overline k_3=0 \Rightarrow \phhi'=0$.
\end{proof}

Every regular $C^2$-curve has a family of Bishop frames and congruent Bishop developments, which can serve as an invariant (up to rotation) of the curve (cf. \cite{bishop}). The existence proof  is omitted as it is not needed for this paper but when a Frenet frame is given, it is straightforward to construct an equivalent Bishop frame.

\begin{theorem}[Bishop Transformation]
\label{bishoptransform}
Let $(\v{T}, \v{N}, \v{B}, \kappa, \tau)$ be a Frenet apparatus of a Frenet curve. Then the family of equivalent Bishop apparatuses \,$(\v T, \v N_1, \v N_2, k_1, k_2)$ is given as follows:
\[ \begin{pmatrix} \v {N}_1\\ \v {N}_2\end{pmatrix}=
 \begin{pmatrix} \cos\phhi& -\sin \phhi  \\ \sin\phhi & \quad \cos \phhi\end{pmatrix}\cdot
\begin{pmatrix} \v {N}\\ \v { B}\end{pmatrix},
\quad %\text{and}\quad
\begin{pmatrix}  {k}_1\\  { k}_2\end{pmatrix}=\kappa  \begin{pmatrix} \cos\phhi\\ \sin \phhi\end{pmatrix},
\quad \phhi=\phhi_0+\int \tau,
\]
with a constant parameter \:$\phhi_0 \in \R$.
This conversion is called {\em Bishop transformation}. 

\end{theorem}
\begin{proof}
The Bishop transformation is an application of theorem \ref{totality} with $k_1=\kappa$, \:$k_2=0$, \:$k_3=\tau$, and $\overline k_3=0$. 
\end{proof}

\begin{theorem}[Inverse Bishop Transformation]
\label{inversebishop}
A regular $C^2$ curve with Bishop apparatus $(\v T, \v N_1, \v N_2, k_1, k_2)$ is a Frenet curve if and only if it has a Bishop development that can be expressed in polar coordinates
\[\begin{pmatrix}  {k}_1\\  { k}_2\end{pmatrix}=\omega  \begin{pmatrix} \cos\phhi\\ \sin \phhi\end{pmatrix}\]
with continuous polar radius $\omega$ and differentiable polar angle $\phhi$. 
A Frenet apparatus can then be constructed as follows: 
$$(\v T, \ \v N=\cos \phhi \v N_1 + \sin \phhi \v N_2,\  \v B=-\sin \phhi \v N_1 + \cos \phhi \v N_2, \  \kappa=\omega, \  \tau=\phhi').$$ 
\end{theorem}
\begin{proof}
Any Frenet curve must have a Bishop development of this form (theorem \ref{bishoptransform}) and if such a Bishop apparatus exists, inverting the Bishop transformation must result in a Frenet apparatus.
% For the conversion from Bishop frame to Frenet frame, theorem \ref{totality} leads to the conditions
% $$ \kappa=k_1\cos\phhi-k_2\sin\phhi, \quad k_2\cos\phhi +k_1\sin\phhi=0, \quad \tau=-\phhi'.$$
% The first two are equivalent to $\kappa=e^{-i\phhi}(k_1+ik_2) \,\Leftrightarrow \,k_1=\kappa\cos\phhi \wedge  k_2=\kappa\sin\phhi$.
\end{proof}

\subsection{Successor Curves}\label{sectionsuccessorcurves}
This section introduces certain transformations of the Frenet apparatus that preserve the Frenet property. In this way, new curves can be derived and solutions to their natural equations be obtained from known cases. 

\begin{defin}[Successor Curve]
Given a unit speed $C^2$ curve $x=x(s)$ with unit tangent $\v T$. A curve $x_1=x_1(s)$ (with the same arc length parameter) that has $\v T$ as principal normal is called a {\em successor curve} of $x$. 
A Frenet frame $ (\v{T_1}, \v{N_1}, \v{B_1})$ is called  {\em successor frame} of the Frenet frame $(\v{T}, \v{N}, \v{B})$ if $\v N_1\equiv \v T$. 
\end{defin}

\begin{theorem}[Successor Transformation of Frenet Apparatus]
\label{successor}
Every Frenet curve has a family of successor curves. Given a Frenet system $\Kal F=(\v{T}, \v{N}, \v{B}, \kappa, \tau)$, the totality of successor systems  $\Kal F_1=(\v{T}_1, \v{N}_1, \v{B}_1, \kappa_1, \tau_1)$ is as follows:

\[
\begin{pmatrix} \v T_1 \\ \v {N_1}\\ \v {B_1}\end{pmatrix}=
 \begin{pmatrix} 0&  -\cos\phhi &  \sin \phhi \\ 1&0&0 \\ 0& \quad\sin\phhi& \cos \phhi  \end{pmatrix}\cdot
\begin{pmatrix} \v T \\ \v {N}\\ \v { B}\end{pmatrix}, \quad
\begin{pmatrix}  \kappa_1\\  \tau_1\end{pmatrix}=\kappa  \begin{pmatrix} \cos\phhi\\ \sin \phhi\end{pmatrix},\quad \phhi(s)\;=\phhi_0+\int_{s_0}^s \tau(\sigma)d\sigma, 
\]
with a constant $\phhi_0\in\R$. The Darboux vector of $\Kal F_1$ is \,$\v D_{\omega 1}=\kappa \v{B}$ and is normalizable with angular speed $\kappa$ (the successor curve's signed Lancret curvature).
\end{theorem}

\begin{proof}
The tangent of a successor curve is a parallel unit vector field relative to the original curve. Any normal component of a Bishop frame thus is a tangent to a successor curve, and vice versa. Therefore, all successor frames can be obtained by applying the Bishop transformation (\ref{bishoptransform}) to any Frenet frame of the original curve and rearranging the resulting Bishop apparatus into a Frenet apparatus (\ref{rearrange} e). The family of successor curves is well-defined in that it does not depend on the choice of Frenet system of the original curve. 
\end{proof}

\begin{theorem}[Predecessor Transformation of Frenet Apparatus]
\label{predecessor}
Given a Frenet apparatus with a Frenet development that can be expressed in polar coordinates, 
\[\Kal F_1= (\v{T_1}, \v{N_1}, \v{B_1},\: \kappa_1=\omega \cos \phhi, \:\tau_1=\omega \sin \phhi)\] ($\omega$ continuous, $\phhi$ differentiable). 
Then 
\:$\Kal F= (\v{T}, \v{N},\v{B},\:\kappa=\omega, \:\tau=\phhi') \:\:\text{with} $
\[ \v{T}=\v N_1,\quad \v{N}=-\cos \phhi \v T_1 + \sin \phhi \v B_1,\quad \v{B}=\v D_1= \sin \phhi \v T_1 + \cos \phhi \v B_1.\]
is a Frenet apparatus, called a {\em predecessor} apparatus of $\Kal F_1$. $\v D_1$ is the unit Darboux vector of $\Kal F_1$ and $\Kal F_1$ is a successor apparatus of $\Kal F$. 

If $\Kal F_1$ is at least $C^2$ and has positive Lancret curvature $\omega_{1+}=\; \sqrt{\kappa_1^2+\tau_1^2}=:\omega$, then it has a uniquely (up to sign) determined, strongly regular predecessor apparatus
\[\v N=\frac{1}{\omega}(- \kappa_1 \v T_1+\tau_1 \v B_1 ), \quad \v B=\frac{1}{\omega}(\tau_1 \v T_1 + \kappa_1 \v B_1), \quad \kappa=\omega>0, \quad \tau=\frac{\kappa_1 \tau_1'-\kappa_1'\tau_1}{\omega^2}.\]
\end{theorem}
\begin{proof}
The predecessor transformation is obviously the inverse of the successor transformation \ref{successor}.
\noindent
If $\kappa_1$ and $\tau_1$ never vanish simultaneously, they can be expressed in polar coordinates with radius $\omega_{1+}$ and angle $\phhi$. 
Observe that $\kappa_1\cos\phhi=\tau_1\sin\phhi$. Differentiating and multiplying with $\omega$ \,yields \,$\kappa_1\tau_1'-\kappa_1'\tau_1=\phhi'\omega^2$ \,and the result for $\tau$ follows. If $\kappa_1>0,\;\phhi=\arctan \tau_1/\kappa_1$. 
\end{proof}
\begin{rem}
The predecessor transformation for the strongly regular case was described in \cite{bilinski1955}. See also \cite{bilinski1963} and \cite{Hoschek1969}.
Note that in general, there is no well-defined predecessor curve for a given Frenet curve because it depends on the choice of principal normal. 
\end{rem}

\begin{rem}[Closed Curves and Periodic Frames]
\label{periodic}
Under what conditions a Frenet curve is closed is an open problem, known as the {\em closed curve problem} (\cite{fenchel1951}).
The existence of a Frenet development that is periodic with respect to the arclength is necessary but not sufficient for the curve to be closed. Note that the arc length parametrization of a closed curve of length $L$ is obviously $L$-periodic but it need not have an $L$-periodic Frenet apparatus, only one that is $2L$-periodic (picture a Frenet frame along a Mobius strip reversing its orientation after each loop). 

The total torsion of a closed curve is a curve invariant (modulo $\pi$) not depending on the choice of the Frenet frame (see \ref{totaltorsion}; \cite{wonglai}: 
15\footnote{Wong and Lai (1967) mistakenly state that the total torsion is a curve invariant modulo $2\pi$.}.)
In the case of a periodic Frenet apparatus,\, $\Kal F(s+P)=\Kal F(s)$ for all $s$, the {\em torsion angle} $\int_{s_0}^{s_0+P}\tau(s)\, ds \mod \pi$ is a curve invariant at least if $\kappa(s_0)\neq 0$. 
It corresponds to the change in the polar angle of the equivalent Bishop development (theorem \ref{inversebishop}). Any Bishop apparatus, and consequently any successor apparatus $\Kal F_1$, is periodic if and only if the total torsion of $\Kal F$ is a rational multiple of $\pi$ (theorem \ref{successor}).

Conversely, if $\Kal F_1$ is periodic with positive Lancret curvature, its Frenet development is a twodimensional closed loop not meeting the origin. Any predecessor apparatus is also periodic and its total torsion corresponds to the polar angle traversed during a full loop (theorem \ref{predecessor}), therefore it is an integer multiple of $2\pi$. If $\kappa_1>0$, the Frenet development cannot loop around the origin and the total torsion of the predecessor curve vanishes (compare \cite{hartl}: 37).
\end{rem}

\begin{rem}[Closed Spherical Curves]
In the case of {\em lines of curvature} (remark \ref{darbouxframe}), the {\em Darboux moving frame} forms a Bishop frame. Because the surface normal in any point is unique up to sign, a closed line of curvature of length $L$ has a $2L$-periodic Bishop frame and its total torsion is $k\pi$ (\cite{fenchel1934}; \cite{menninger}: 61; \cite{qinli}). %Qin and Li 2002
Spherical curves are a special case of lines of curvature. Their Bishop development lies on a straight line not meeting the origin (\cite{bishop}). A closed spherical Frenet curve of length $L$ has $L$-periodic Frenet and Bishop frames and, as is well known, its total torsion vanishes (\cite{menninger}: 61-64).
\end{rem}

\newpage
\section{Plane Curves, Helices and Slant Helices}
\label{sectionapplied}

Wo now apply the theory of successor curves to solve the natural equations of the helix and slant helix.
\subsection{Plane Curves}

\begin{theorem}[Plane Curves]

\label{plane}

A regular $C^2$ space curve lies on a plane if and only if it has a generic Frenet apparatus with constant binormal and vanishing torsion.

\noindent Given continuous functions $\kappa=\kappa(s)$ and $\tau\equiv0$ defined on an interval $s_0 \in I \subset \R$. Let $\quad{\mit\Omega}  (s)=\int_{s_0}^s\kappa(\sigma)d\sigma$ and 
\[\v{T}_P(s)=\begin{pmatrix} \cos{\mit\Omega(s)}\\ \sin{\mit\Omega(s)} \\0 \end{pmatrix}, \quad
\v{N}_P(s)=\begin{pmatrix}   -\sin{\mit\Omega(s)}\\  \quad\cos{\mit\Omega(s)} \\0 \end{pmatrix},\quad
\v{B}_P=\begin{pmatrix} 0\\ 0 \\1 \end{pmatrix}.\]
Then $\Kal F_P= (\v{T}_P, \v{N}_P, \v{B}_P, \kappa, \tau)$ is a generic Frenet apparatus of the plane curve characterized by  $\kappa$ and $\tau$. Its tangent and principal normal images lie on the equator, its binormal is constant and its axis of motion is fixed, with unit Darboux vector $\v D_P=\v{B}_P$.
\end{theorem}

\begin{proof}
Let $\v{x}(s)$ be a regular plane curve and $\v{B_0}$ a unit normal of the plane. If $\v{x}(s_0)$ is a curve point, $\v{x}(s)-\v{x}(s_0) \perp \v B_0 \Rightarrow \inner{\v{x}'(s)}{\v{B}_0}=\inner{\v{T}(s)}{\v{B}_0}=0 \Rightarrow \inner{\v{T}'(s)}{\v{B}_0}=0$. By definition \ref{frenetcurves}, $\v B_0$ is a binormal of the curve and the curve is a Frenet curve. With a constant binormal, $\tau=\inner{\v{B}_0'}{\v{N}} \equiv 0 $.

Finally, it is easy to verify that $\Kal F_P$ satisfies the definition of a generic Frenet apparatus and determines a plane curve. Its Darboux vector $\v D_{\omega P}=\kappa \v B_P$ is normalizable.
\end{proof}

\subsection{General Helices}

\begin{defin}[General Helix, Slope Line]
\label{helixdefin}
A regular $C^2$-curve is called a {\em helix} or {\em general helix} (also {\em slope line} or {\em curve of constant slope}) if at every point its tangent vector makes a constant angle (the {\em slope angle}) $\theta \in \; (0, \pi/2)$ with a fixed direction, represented by a unit vector $\v D_0$ called the {\em slope vector}.
\end{defin}
\begin{rem}
Straight lines ($\theta=0$) and plane curves ($\theta=\pi/2$) are excluded.
\end{rem}
\begin{theorem}[General Helices]
\label{helix}
A curve is a general helix if and only if it is a successor curve of a plane curve, and not itself a plane curve.
\end{theorem}
\begin{proof}
Given a helix with unit tangent $\v T$ and slope vector $\v D_0$. It is $\inner{\v T}{\v D_0}=const. \Leftrightarrow \inner{\v T'}{\v D_0}=0$. 
$\v N=\csc\theta\, \v D_0 \times \v T$ is a unit normal vector that normalizes $\v T'$ and therefore is a principal normal. $\int\v N$ is a plane curve with plane normal $\v D_0$ and tangent $\v N$ and by definition, the helix is its successor curve.

Conversely, let $\v T$ be the tangent of a plane curve with unit plane normal $\v B_0$ and $\v T_1$ the tangent of a successor curve. $\inner{\v T_1'}{\v B_0}=\inner{\kappa_1 \v T}{\v B_0}=0 \Leftrightarrow \inner{\v T_1}{\v B_0}=const.$
\end{proof}

\begin{theorem}[Frenet Apparatus of General Helix]
\label{helix2}
A regular $C^2$ space curve is a general helix with slope angle  $\theta \in \; (0, \pi/2)$ if and only if it has a Frenet development satisfying $$\tau=\cot \theta \kappa.$$ 

Given a constant $\theta \in \, (0, \pi/2)$ and continuous functions $\kappa$, ${\kappa_H=\sin\theta\,\kappa}$, $\tau_H=\cot \theta \, \kappa_H=\cos\theta\, \kappa: I \mapsto \R, \ s_0 \in I$. Let ${{\mit\Omega}  (s)=\int_{s_0}^s\kappa(\sigma)d\sigma}$ and 
\[\v{T}_H(s)=\begin{pmatrix} \quad\sin\theta \sin{\mit\Omega(s)}\\ -\sin\theta\cos{\mit\Omega(s)} \\ \cos\theta \end{pmatrix}, \;
\v{N}_H(s)=\begin{pmatrix}  \cos{\mit\Omega(s)}\\ \sin{\mit\Omega(s)} \\0  \end{pmatrix},\;
\v{B}_H(s)=\begin{pmatrix} -\cos \theta \sin{\mit\Omega(s)}\\ \quad\cos \theta \cos{\mit\Omega(s)}\\ \sin\theta \end{pmatrix}.\]

Then $\Kal F_H= (\v{T}_H, \v{N}_H, \v{B}_H, \kappa_H, \tau_H)$ is a generic Frenet apparatus of a general helix with slope angle $\theta$. Its tangent has constant slope $\theta$, its binormal has constant slope $\pi/2-\theta$ and its principal normal image lies on a great circle. Its axis of motion is fixed in the direction of the slope vector $\v D_0$, with unit Darboux vector $\v D_H=\v D_0$ and angular speed $\kappa$. 
\end{theorem}
\begin{proof}  
Any helix is successor curve of a plane curve and has a Frenet apparatus $\Kal F_H$ that is successor apparatus (theorem \ref{successor}) of the generic Frenet apparatus of a plane curve (theorem \ref{plane}). Every choice of $\phhi_0=\theta\in \; (0, \pi/2)$ results in a distinct successor curve. $\Kal F_H$ is generic as $\tau_H$ vanishes whenever $\kappa_H$ vanishes.
\end{proof}

\subsection{Slant Helices}
\begin{defin}[Slant Helix]
A Frenet curve is called a {\em slant helix} (\cite{izumya}) if it has a principal normal that has constant slope $\theta \in \; (0, \pi/2)$. 
\end{defin}

Slant helices by definition are exactly the successor curves of proper general helices and a Frenet frame is obtained by applying the successor transformation. A parametrization of the unit tangent is now given. For convenience we'll denote $m:=\cot \theta, \,n:=m/\sqrt{1+m^2}=\cos\theta, \,n/m=\sin\theta$. 

\begin{theorem}[Frenet Apparatus of Slant Helix]\label{slanthelix}
A regular $C^2$ space curve is a slant helix if and only if it has a Frenet development satisfying 
\[\kappa_{SH}=\frac{1}{m} \,\phhi' \cos\phi,\;\quad\tau_{SH}=\frac{1}{m} \,\phhi' \sin\phi\] 
with a differentiable function $\phi(s)$ and $m=\cot\theta\neq 0$.

Given such a Frenet development and let 
$\mit\Omega(s):=\phi(s)/n,\, \lambda_1:=1-n, \lambda_2:= 1+n$ with $n=\cos\theta$. Then the tangent vector of the slant helix thus characterized can be parametrized as follows:
\[\v{T}_{SH}(s)=\frac{1}{2} \begin{pmatrix}
 \lambda_1\cos\lambda_2{\mit\Omega(s)}  +\lambda_2\cos\lambda_1{\mit\Omega(s)}  \\
 \lambda_1\sin\lambda_2{\mit\Omega(s)}  +\lambda_2\sin\lambda_1{\mit\Omega(s)}  \\
    2\frac{n}{m}\sin n{\mit\Omega(s)}   \end{pmatrix}\]

The slant helix is the successor curve of a helix with curvature \,$\kappa_H=\phi'/m$ and torsion \,$\tau_H=\phi'$.
The regular arcs traced by $\v{T}_{SH}$ are spherical helices and the slant helix has a Darboux vector with constant slope $\pi/2-\theta$ and speed of precession $\omega=\kappa_H$. 
Slant helices are generic.
\end{theorem}

\begin{proof}
The successor transformation  (theorem \ref{successor}) of $\Kal F_H$ gives \,$\phi=\phi_0+m\int\kappa_H$ and \,$\kappa_{SH}=\kappa_H\cos\phi=m^{-1}\phi'\cos\phi$, $\tau_{SH}=m^{-1}\phi'\sin\phi$. Further, $\mit\Omega=\mit\Omega_0+ m/n\int\kappa_H$. The integration constant $\mit\Omega_0$ (omitted in theorem \ref{helix2}) can be set arbitrarily because it only effects a rotation of the frame of the helix. We set $\mit\Omega_0=\phi_0/n$ to get $\phi=n\,\mit\Omega$. 
For the successor tangent, we have $\v{T_{SH}}=-\cos\phhi \v{N}_H + \sin \phhi \v{\v{B}}_H$. Therefore
\[\v{T}_{SH}= \begin{pmatrix} -\cos n\mit\Omega \cos \mit\Omega - n\sin n\mit\Omega\sin\mit\Omega\\
 -\cos n\mit\Omega \sin \mit\Omega + n\sin n \mit\Omega\cos\mit\Omega\\
\frac{n}{m}\sin n \mit\Omega \end{pmatrix}\]
which can be simplified as above (the first two components have been reflected for convenience). Parametrizations of principal normal and binormal can be obtained in the same way but are omitted here. 

The regular segments of the tangent image are slope lines with the principal normal as tangent.
By construction, the unit Darboux vector of the slant helix is the binormal of a helix, which has constant slope (cf. theorem \ref{successor}). The Darboux vector is $\kappa_H B_H$ and the function $\omega=\kappa_H$ is the signed Lancret curvature. The Frenet apparatus thus constructed is generic, since along line segments, $\phhi'=0$ and therefore curvature and torsion both vanish.
\end{proof}

\begin{rems}[Natural Equations of Slant Helix]~\\ 
\label{slanthelixrems}
\indent 1. The natural equations for slant helices derived above are of unrestricted generality. Any choice of a continuous function $\phi(s)$ and constant $\theta\in(0, \pi/2)$ gives rise to a well-defined slant helix.

2. It is easy to verify that under the assumption $\kappa \neq 0$, the natural equations of the slant helix are equivalent to the condition
\[ \frac{\kappa^2}{(\kappa^2+\tau^2)^{3/2}} \left( \frac{\tau}{\kappa}\right)' = m.\]
This expression is equivalent to the {\em geodesic curvature} of the spherical image of the principal normal (see \cite{izumya}), which is part of a circle.

3. Solving the differential equation \, $\phhi'\cos\phhi=m\kappa$ for $\phhi$ yields \, $\cos\phhi \,d\phhi=m\,\kappa \,ds\;\Rightarrow \;\sin\phhi(s)=m\int\kappa\, ds$. Similarly, \,$\cos \phi(s)=-m\int\tau ds$. Thus the plane curve with total curvature and total torsion as coordinates is circular. Further, if $\phi \mod 2\pi$ is a periodic function (and consequently the Frenet apparatus is periodic), then total curvature and total torsion over each period vanish.

4. Restricting \,$\phhi$ to $(-\pi/2, \pi/2)$, the torsion can be expressed as $\tau=\kappa\tan\phhi=\kappa\frac{\sin\phhi}{\sqrt{1-\sin^2\phhi}}$. It follows (cf. \cite{ali}):
\[\tau(s)=\,\kappa(s)\frac{m\int\kappa(s)ds}{\sqrt{1-m^2(\int\kappa(s)ds)^2}}.\]
Given an arbitrary curvature $\kappa$, a torsion function can be constructed to create a slant helix, provided that the domain of $s$ is appropriately restricted. 

5. An interesting application is the construction of a {\em Salkowski curve} (\cite{salkowski1909}, \cite{monterde}, \cite{ali}), a slant helix with constant curvature:
\[\kappa_S(s)\equiv 1, \quad \tau_S(s) = \frac{m s}{\sqrt{1-m^2 s^2}}, \quad m\neq 0, \quad  s\in\, (-1/|m|, +1/|m|).\]
It follows $\,\sin\phi=ms, \,\phi=\arcsin ms, \,\cos\phi=\sqrt{1-m^2}$.
The Salkowski curve is the successor curve of its unit tangent vector, which is a spherical helix and has the same arc length. Its curvatures are
\[\kappa_H=\frac{1}{\sqrt{1-m^2s^2}}, \quad\tau_H=m\,\kappa_H.\]

6. The case \,$\phi'=const.$ gives rise to the class of {\em curves of constant precession} (\cite{scofield}). Setting \,$\phi(s) :=\mu s$ and \,$\omega:=\mu/m$ and \,$\alpha:=\mu/n=\sqrt{\omega^2+\mu^2}$, we substitute $\mit\Omega=\alpha s$ and $\lambda_1=(\alpha-\mu)/\alpha$,\, $\lambda_2=(\alpha+\mu)/\alpha$ in the formulae for the slant helix.
\end{rems}

%\\pagebreak
\begin{theorem}[Curves of Constant Precession]
 A {\em curve of constant precession} has curvature, torsion and tangent parametrization as follows:
\[\kappa_{CP}=\omega\cos\mu s, \; \tau_{CP} = \omega \sin \mu s \]
\[\v{T}_{CP}= \frac{1}{2\alpha} \begin{pmatrix}
 (\alpha - \mu)\cos(\alpha+\mu)s + (\alpha + \mu)\cos(\alpha-\mu)s    \\
 (\alpha - \mu)\sin(\alpha+\mu)s +  (\alpha + \mu)\sin(\alpha-\mu)s  \\
    2\omega\sin\mu s  \end{pmatrix}\]
\[\text{with arbitrary constants} \;\omega, \mu \neq 0 \; \text{and }\alpha=\sqrt{\omega^2+\mu^2}.\]

\noindent
An explicit arclength parametrization of the curve is obtained by elementary integration. The resulting curve of constant precession lies on a hyperboloid of one sheet, and it is closed if and only if 
$\cos\theta=\mu/\alpha$
is rational (figure 2 in \cite{scofield}). 
\end{theorem}

\begin{rem}
Curves of constant precession are exactly the successor curves of circular helices with curvature $\omega$ and torsion $\mu$, which in turn are successor curves of circles. 
Closed plane curves have zero total torsion and their successor helices have periodic Frenet frames (\ref{periodic}) but are not closed. The Frenet frame of the circular helix has period $L=2\pi/\alpha$ and its total torsion is $2\pi\mu/\alpha$ and so the successor apparatus is periodic if and only if $\mu/\alpha$ is rational. In this case, we obtain a closed curve. Both the total curvature and total torsion of a closed curve of constant precession vanish (consistent with remark \ref{slanthelixrems}.4). 

Total curvature has a geometric interpretation as follows. The absolute total curvature on a strongly regular curve segment is equivalent to the length of the arc traced by the tangent on the sphere. 
Vanishing total curvature implies that the arc lengths of the tangent image between the cusps (inflection points) add up to zero, when each is assigned the sign of the curvature. In the case of curves of constant precession, each "upward" arc of the tangent image is balanced out by a "downward" arc of the same length (figure 1 in \cite{scofield}). 
\end{rem}

\subsection*{Outlook}
In this paper, a successor transformation was introduced for Frenet curves and it was shown that the successor curves of plane curves are helices, and their successor curves slant helices. This procedure could be repeated indefinitely to give rise to new, yet unexplored classes of curves. In particular, applying the transformation to a closed curve of constant precession should result in a curve with periodic Frenet apparatus and potentially also closed.

%\pagebreak
\bibliographystyle{apalike}
\bibliography{curves}

\begin{thebibliography}{}

\bibitem[Ali, 2012]{ali}
Ali, A.~T. (2012).
\newblock Position vectors of slant helices in {E}uclidean 3-space.
\newblock {\em Journal of the Egyptian Mathematical Society}, 20:1--6, DOI: \href{http://www.sciencedirect.com/science/article/pii/S1110256X11000289}{10.1016/j.joems.2011.12.005}.

\bibitem[Bilinski, 1955]{bilinski1955}
Bilinski, S. (1955).
\newblock Eine {V}erallgemeinerung der {F}ormeln von {F}renet und eine {I}somorphie gewisser 
{T}eile der {D}ifferentialgeometrie der {R}aumkurven.
\newblock {\em Glasnik Math.-Fiz. i Astr.}, 10:175--180.

\bibitem[Bilinski, 1963]{bilinski1963}
Bilinski, S. (1963).
\newblock \"{U}ber eine {E}rweiterungsm\"{o}glichkeit der {K}urventheorie.
\newblock {\em Monatshefte für Mathematik}, 67:289--302, \href{http://eudml.org/doc/177219}{eudml.org/doc/177219}.

\bibitem[Bishop, 1975]{bishop}
Bishop, R.~L. (1975).
\newblock There is more than one way to frame a curve.
\newblock {\em Amer. Math. Monthly}, 82:246--251, \href{http://www.jstor.org/stable/2319846}{www.jstor.org/stable/2319846}.

\bibitem[Camci et~al., 2013]{Camci2013}
Camci, C., Kula, L., and Altinok, M. (2013).
\newblock On spherical slant helices in euclidean 3-space.
\newblock {\em arXiv:1308.5532 [math.DG]}.

\bibitem[Fenchel, 1934]{fenchel1934}
Fenchel, W. (1934).
\newblock \"{U}ber einen {J}acobischen {S}atz der {K}urventheorie.
\newblock {\em T\^{o}hoku Math. J.}, 39:95--97.

\bibitem[Fenchel, 1951]{fenchel1951}
Fenchel, W. (1951).
\newblock On the differential geometry of closed space curves.
\newblock {\em Bull. Amer. Math. Soc.}, 57:44--54, \href{http://projecteuclid.org/euclid.bams/1183515801}{http://projecteuclid.org/euclid.bams/1183515801}.


\bibitem[Hartl, 1993]{hartl}
Hartl, J. (1993).
\newblock Zerlegung der {D}arboux-{D}rehung in zwei ebene {D}rehungen.
\newblock {\em Journal of Geometry}, 47:32--38, \href{http://link.springer.com/article/10.1007\%2FBF01223802}{link.springer.com/article/10.1007\%2FBF01223802}.

\bibitem[Hoppe, 1862]{Hoppe1862}
Hoppe, R. (1862).
\newblock \"{U}ber die {D}arstellung der {C}urven durch {K}rümmung und {T}orsion.
\newblock {\em Journal für die reine und angewandte Mathematik}, 60:182--187, \href{http://eudml.org/doc/147848}{http://eudml.org/doc/147848}.

\bibitem[Hoschek, 1969]{Hoschek1969}
Hoschek, J. (1969).
\newblock Eine {V}erallgemeinerung der {B}\"oschungsfl\"achen.
\newblock {\em Mathematische Annalen}, 179:275--284, \href{http://eudml.org/doc/161782}{http://eudml.org/doc/161782}.

\bibitem[Izumiya and Takeuchi, 2004]{izumya}
Izumiya, S. and Takeuchi, N. (2004).
\newblock New {S}pecial {C}urves and {D}evelopable {S}urfaces.
\newblock {\em Turk J Math}, 28:153--163,
\href {http://journals.tubitak.gov.tr/math/issues/mat-04-28-2/mat-28-2-6-0301-4.pdf}{journals.tubitak.gov.tr/math/issues/mat-04-28-2/mat-28-2-6-0301-4.pdf}.

\bibitem[K\"uhnel, 2006]{kuehnel}
K\"uhnel, W. (2006).
\newblock {\em Differential Geometry: Curves - Surfaces - Manifolds, Second Edition}.
\newblock American Mathematical Society.

\bibitem[{Kula} et~al., 2010]{kula2010}
{Kula}, L., {Ekmekci}, N., {Yayl{\i}}, Y., and {\.Ilarslan}, K. (2010).
\newblock {Characterizations of slant helices in Euclidean 3-space.}
\newblock {\em {Turk. J. Math.}}, 34(2):261--274,
DOI: \href{http://dx.doi.org/10.3906/mat-0809-17}
{10.3906/mat-0809-17}.

\bibitem[{Kula} and {Yayli}, 2005]{kula2005}
{Kula}, L. and {Yayli}, Y. (2005).
\newblock {On slant helix and its spherical indicatrix.}
\newblock {\em {Appl. Math. Comput.}}, 169(1):600--607,
DOI: \href{http://dx.doi.org/10.1016/10.1016/j.amc.2004.09.078}
{10.1016/j.amc.2004.09.078}.

\bibitem[Menninger, 1996]{menninger}
Menninger, T. (1996).
\newblock \"{U}ber die {D}arstellung von {R}aumkurven durch ihre {I}nvarianten, \href{http://www.slideshare.net/amenning/uber-die-darstellung-von-raumkurven-durch-ihre-invarianten}{www.slideshare.net/amenning/uber-die-darstellung-von-raumkurven-durch-ihre-invarianten}.

\bibitem[Monterde, 2009]{monterde}
Monterde, J. (2009).
\newblock Salkowski curves revisted: A family of curves with constant curvature
and non-constant torsion.
\newblock {\em Computer Aided Geometric Design}, 26:271–278, DOI: \href{http://dx.doi.org/10.1016/j.cagd.2008.10.002}
{10.1016/j.cagd.2008.10.002}

\bibitem[Nomizu, 1959]{nomizu}
Nomizu, K. (1959).
\newblock On {F}renet equations for curves of class ${C}^{\infty}$.
\newblock {\em T\^ohoku Math. J.}, 11:106--112, \href{http://projecteuclid.org/euclid.tmj/1178244631}{projecteuclid.org/euclid.tmj/1178244631}.

\bibitem[Qin and Li, 2002]{qinli}
Qin, Y.-A. and Li, S.-J. (2002).
\newblock Total {T}orsion of {C}losed {L}ines of {C}urvature.
\newblock {\em Bull. Austral. Math. Soc.}, 65:73--78,
DOI: \href{http://dx.doi.org/10.1017/S0004972700020074}
{10.1017/S0004972700020074}

\bibitem[Salkowski, 1909]{salkowski1909}
Salkowski, E. (1909).
\newblock Zur {T}ransformation von {R}aumkurven. 
\newblock {\em Mathematische Annalen} 
66: 517--557,
\href{http://eudml.org/doc/158392}{eudml.org/doc/158392}.

\bibitem[Scofield, 1995]{scofield}
Scofield, P.~D. (1995).
\newblock Curves of {C}onstant {P}recession.
\newblock {\em Amer. Math. Monthly}, 102:531--537, \href{http://www.jstor.org/stable/2974768}{www.jstor.org/stable/2974768}.

\bibitem[Wintner, 1956]{wintner}
Wintner, A. (1956).
\newblock On {F}renet's {E}quations.
\newblock {\em Amer. J. Math.}, 78:349--355, \href{http://www.jstor.org/stable/2372520}{www.jstor.org/stable/2372520}.

\bibitem[Wong and Lai, 1967]{wonglai}
Wong, Y.-C. and Lai, H.-F. (1967).
\newblock A {C}ritical {E}xamination of the {T}heory of {C}urves in {T}hree
  {D}imensional {D}ifferential {G}eometry.
\newblock {\em T\^ohoku Math. J.}, 19:1--31, \href{http://projecteuclid.org/euclid.tmj/1178243344}{projecteuclid.org/euclid.tmj/1178243344}.

\end{thebibliography}

\end{document}